\documentclass[10pt]{article}
\usepackage{amsfonts}
\begin{document}
\newtheorem{theorem}{Theorem}[section]
\newtheorem{lemma}[theorem]{Lemma}
\newtheorem{result}[theorem]{Result}
\newcommand{\ord}{\ensuremath{\mathrm{ord}}}
\newcommand{\proofbox}{\hspace*{\fill}\quad$\Box$}
\newcommand{\PLR}{primitive $\lambda$-root\ }
\newcommand{\PLRs}{primitive $\lambda$-roots\ }
\newcommand{\mod}[1]{\ensuremath{(\mathrm{mod}\,#1)}}

\newtheorem{unnumber}{}
\renewcommand{\theunnumber}{\relax}
\newtheorem{prepf}[unnumber]{Proof}
\newenvironment{pf}{\prepf\rm}{\proofbox\endprepf}

\newcommand{\zn}{\ensuremath{\mathbb{Z}_n}}
\newcommand{\zz}[1]{\mathbb{Z}_{#1}}
\newcommand{\un}{\ensuremath{\mathbb{U}_n}}
\newcommand{\uu}[1]{\mathbb{U}_{#1}}
\newcommand{\bbZ}{\ensuremath{\mathbb{Z}}}
\newcommand{\cy}[2]{\langle#1\rangle_{#2}}
\newcommand{\GF}{\mathop{\mathrm{GF}}}

\title{Three-factor decompositions of $\un$ with\\
        the three generators in arithmetic progression}
\author{P.\ J.\ CAMERON$^\textup{\small a}$ and 
        D.\ A.\ PREECE$^\textup{\small a,b}$ 
        \bigskip\\
        $^\textup{\small a}$\textit{School of Mathematical Sciences,}\\
        \textit{Queen Mary University of London,}\\
        \textit{Mile End Road, London  E1 4NS, UK}\\
        \texttt{P.J.Cameron@qmul.ac.uk, D.A.Preece@qmul.ac.uk}
        \bigskip\\
        $^\textup{\small b}$\textit{Institute of Mathematics,
                       Statistics and Actuarial Science,}\\
        \textit{Cornwallis Building, University of Kent,}\\
        \textit{Canterbury, Kent  CT2 7NF, UK}
        }
\maketitle
\begin{abstract}
Irrespective of whether $n$ is prime, prime power with exponent $> 1$, 
or composite, the group $\un$ of units of $\zn$ can sometimes be 
obtained as $\un = \langle x \rangle \times \langle x+k \rangle$ 
$\times \langle x+2k \rangle$ where $x,k \in \zn$.
Indeed, for many \mbox{values} of $n$, many distinct $3$-factor decompositions
of this type exist.   The circum\-stances in which such decompositions
exist are examined. Many \mbox{decompositions} have additional interesting
properties. We also look briefly at decompositions of the multiplicative
groups of finite fields.\\

 
\end{abstract}

\section{Introduction}
\label{s1}

An element of $\zn$ is a \emph{unit} of $\zn$ if $x$ and $n$ are
co-prime.   If the prime-power decomposition of $n$ is
$n = p^\alpha q^\beta r^\gamma \cdots$ where $p, q, r, \ldots\ $ are 
distinct primes, then the number of units is given by Euler's
\emph{totient function} $\phi_n = |\un| =$ 
$(p-1)p^{\alpha - 1} \cdot (q-1)q^{\beta - 1} 
                     \cdot (r-1)r^{\gamma - 1} \cdots\,$
where $\un$ denotes the group of units \cite[Chap.~5]{JJ} 
We recall the structure of $\un$:
\begin{itemize}
\item If $n=p^\alpha q^\beta r^\gamma\cdots$, where $p,q,r,\ldots$ are primes,
then 
\[\un=\mathbb{U}_{p^\alpha}\times \mathbb{U}_{q^\beta}\times
\mathbb{U}_{r^\gamma}\times\cdots\ ;\]
\item If $p$ is an odd prime, then $\mathbb{U}_{p^\alpha}$ is cyclic of
order $\phi_{p^\alpha}=p^{\alpha-1}(p-1)$;
\item $\mathbb{U}_{2^\alpha}$ is cyclic of order $2^{\alpha-1}$ if
$\alpha\le2$, and is isomorphic to $C_2\times C_{2^{\alpha-2}}$ otherwise.
\end{itemize}
In particular, although $\un$ may be expressible as a direct product of
cyclic groups in many different ways, the smallest number of factors is
equal to the number $k$ of prime divisors of $n$ if $n$ is odd, or $k-1$
if $n$ is twice odd, or $k$ if $n$ is four times odd, or $k+1$ otherwise.
The largest number of factors depends on the prime decompositions of $p-1$
for the prime divisors $p$ of $n$. In particular, if $n$ is prime, then
the maximum number of factors is the number of distinct prime divisors of
$n-1$.

We are mainly concerned here with cases where $n$ is odd and $\un$ is the
product of three cyclic factors. As noted, this forces strong conditions on
$n$: it should have  at most three prime divisors; and if $n$ is prime, then
$n-1$ should have at least three prime divisors.

If we can write 
$\un = \langle x \rangle \times \langle y \rangle \times \langle z \rangle$
where $\langle x \rangle$ denotes the subgroup generated by $x$ and the
symbol $\,\times\,$ connotes a direct product, then we have a
\emph{three-factor decomposition} of $\un$.   If the orders \mod{n}
of $x$, $y$ and $z$ are respectively $a$,~$b$~and~$c$, and we need to
specify them succinctly, we use the notation
\[ \un = \cy{x}{a}\times\cy{y}{b}\times\cy{z}{c} \]
and $\ord_n(x) = a$ \textit{etc.}

Irrespective of whether $n$ is prime, a prime power with exponent ${}> 1$,
or a composite, we can sometimes write
\begin{eqnarray}
\un = \langle x \rangle \times \langle x+k \rangle 
                         \times \langle x+2k \rangle
\label{sct1}
\end{eqnarray}
for some $x,k\in\zn$, so that the generators are in
arithmetic progression (AP).  We then have a \emph{three-factor} AP
\emph{decomposition} of $\un$, which we abbreviate to a 
``3AP decomposition'' of $\un$.   Notable examples~are
\[ \uu{61} = \cy{9}{5} \times \cy{11}{4} \times \cy{13}{3}\ \,        \]
and
\[ \uu{911} = \cy{196}{13} \times \cy{550}{10} \times \cy{904}{7}\ ,  \]
where the orders, like the generators, are in arithmetic progression;
\[ \uu{455} = \cy{92}{4} \times \cy{93}{12} \times \cy{94}{6}\ ,      \]
where the generators are consecutive integers and the orders are all even,
each being one less than a prime factor of $455$;
\[ \uu{91}  = \cy{9}{3} \times \cy{18}{12} \times \cy{27}{2}
            = \cy{87}{6} \times \cy{83}{4} \times \cy{79}{3}          \]
and
\[ \uu{65}  = \cy{61}{3} \times \cy{57}{4} \times \cy{53}{4}\ ,       \]
where $x = k$ in each of the 3AP decompositions; and
\[ \uu{703} = \cy{700}{9} \times \cy{701}{36} \times \cy{702}{2}\ ,   \]
where the generators are respectively ${-3}$, ${-2}$ and ${-1}$
\mod{703 = 19 \times 37}.

Our specification (1) may be realised even if the values $x$, $x+k$ and
$x+2k$, when reduced \mod{n} to lie in the interval $[1,n-1]$, are
not in arithmetic progression in $\bbZ$.   Thus we have
\[ \uu{31} = \cy{30}{2} \times \cy{2}{5} \times \cy{5}{3}   \]
with $k \equiv 3$ \mod{31}.   We could equally have written this as
\[ \uu{31} = \cy{5}{3} \times \cy{2}{5} \times \cy{30}{2}   \]
with $k \equiv 28 \equiv {-3}$ \mod{31}.   Faced with such a choice
between two equivalent representations, we leave ourselves free to choose
whichever seems the more convenient in the context in which it arises.
A rule to make the choice with $0 < k < (n-1)/2$ would be unsatisfactory,
especially as some decompositions~(\ref{sct1}) fall into infinite series 
within which $k$ lies variously in $(0,(n-1)/2)$ and in $((n-1)/2,n)$.

Examples where the two \textbf{outer} generators differ by $1$ include
\[ \uu{275} = \cy{136}{5} \times \cy{-1}{2} \times \cy{137}{20}      \]
and
\[ \uu{775} = \cy{386}{15} \times \cy{-1}{2} \times \cy{387}{20}\ .  \]

We now note three possibilities:
\begin{description}
\item{[A]}
For some values of $n$, different 3AP decompositions of $\un$ may 
arise for different factorisations $\phi_n = a \cdot b \cdot c$.   Thus we have
\begin{eqnarray*}
\uu{211} & = &
\cy{15}{6} \times \cy{107}{5} \times \cy{199}{7} \\
                 & = &
\cy{58}{7} \times \cy{134}{15} \times \cy{210}{2} \\
                 & = &
\cy{196}{3} \times \cy{203}{35}\times\cy{210}{2} \ .
\end{eqnarray*}
\item{[B]} 
For a fixed factorisation $\phi_n = a \cdot b \cdot c$ for a fixed $n$,
we may have different 3AP decompositions 
$\langle x \rangle \times \langle y \rangle \times \langle z \rangle$
of $\un$ where the values of $\ord_n(y)$ are different members of
$\{a, b, c\}$.   Thus we have
\begin{eqnarray*}
\uu{31} & = &
\cy{30}{2} \times \cy{2}{5} \times \cy{5}{3} \\
                & = &
\cy{25}{3} \times \cy{30}{2} \times \cy{4}{5}   
\end{eqnarray*}
and
\begin{eqnarray*}
\uu{547} & = &
\cy{40}{3} \times \cy{172}{26} \times \cy{304}{7} \\
                 & = &
\cy{40}{3} \times \cy{544}{7} \times \cy{501}{26} \\
                 & = &
\cy{520}{7} \times \cy{40}{3} \times \cy{107}{26}\ .
\end{eqnarray*}
\item{[C]}
We may have different 3AP decompositions of $\un$ for a fixed
ordering of the terms of a fixed factorisation $\phi_n = a \cdot b \cdot c$
for a fixed~$n$.   Thus we~have
\begin{eqnarray*}
\uu{191} & = &
\cy{39}{5} \times \cy{190}{2} \times \cy{150}{19} \\
                 & = &
\cy{184}{5} \times \cy{190}{2} \times \cy{5}{19}
\end{eqnarray*}
where $184 \equiv 39^2$ \mod{191} and $5 \equiv 150^{-2}$ \mod{191}.
\end{description}

\paragraph{Problem 1} Find a series of primes that behave like~$191$,
with $2$ as the order of the middle generator.   (Contenders for inclusion
in the series are $n = 191,\ 271$ and $523$.   A possible series with
the order $2$ for an outer generator might cover $n = 331,\ 379,\ 443$
and $647$.)

\medskip

Clearly, a 3AP decomposition of $\un$, where $n$ is prime and
$n > 4$, cannot exist if the
prime-power decomposition of~$n-1$ contains fewer than 3 distinct primes.

Sufficient conditions for the existence of 3AP decompositions of $\un$ 
seem to be elusive.   Thus only computer search has established that,
in the range $n < 300$, a 3AP decomposition of $\un$ does not exist for 
any of the values $n = 71,\ 127,\ 139,\ 223$ and~$277$.   (We here
exclude the ``weak'' 3AP decompositions defined in \S\ref{s2} below.)

For any $n$ with $n > 4$, there is a \emph{primitive root} of $n$ (an
element from $\un$ that generates all members of $\un$)
if and only if $n$ is an odd prime power or twice
an odd prime power.   In general we write $\lambda_n$ for the maximum 
order of a member of $\un$; if $n$ is odd, with prime power decomposition
$n = p^\alpha q^\beta r^\gamma \cdots$, then
\[ \lambda_n = \textup{lcm}((p-1)p^{\alpha - 1},\ (q-1)q^{\beta - 1},
                          \ (r-1)r^{\gamma - 1},\ \ldots)\ .           \] 
We write $\xi_n = \phi_n / \lambda_n$;
as shown in \cite[\S6]{cap}, $\xi_n$ is even if greater than~$1$.

\paragraph{Problem 2} Is there an upper bound on the number of $3$AP
decompositions of $\uu{n}$ in terms of $\xi(n)$? Conversely, for a given
value of $\xi(n)=m$, is it always possible to find $n$ with no $3$AP
decompositions?

Empirically we have found a tendency for larger values of $\xi(n)$ to be
associated with larger numbers of decompositions. The table below,
obtained by computer, gives $D$,
the maximum number of 3AP decompositions of $\uu{n}$, where $n\le1000$ and
$\xi(n)$ is prescribed.

\[\begin{array}{c|rrrrrrrrrrrr}
\xi(n) & 1 & 2 & 4 & 6 & 8 & 10 & 12 & 16 & 18 & 20 & 24 & 36\\\hline
\phantom{|^|}D & 10 & 18 & 96 & 182 & 288 & 262 & 496 &
384 & 276 & 204 & 540 & 2088
\end{array}\]

\section{$n$ prime}
\label{s2}

For $n$ prime, the multiplicative group $\un$ is cyclic, and so if it is
expressed as a direct product, the factors must have pairwise co-prime orders.

\subsection{The case $n-1 = 2 \cdot 3 \cdot m$}
\label{ss21}

We first prove three theorems that apply for prime values $n$ such
that the factors in a 3AP decomposition of $\un$ have orders
$2$, $3$ and $m$, where $2$, $3$ and $m$ are pairwise co-prime.
The first of these theorems is closely linked to Theorem~2.7 of
\cite{consec}. We begin with some preliminary remarks.

Our assumption on $n$ implies that $n\equiv 7$ or $31$ \mod{36},
and $n>7$. In particular, since $n\equiv3$ \mod{4}, the
quadratic residues have odd order, and the non-residues have even order.
In our Theorems, we will be interested in the solutions of the quadratic 
equation $x^2+3x+3=0$ in $\zn$. Its discriminant is $-3$, which 
(by Quadratic Reciprocity \cite[\S7.4]{JJ}) is a square in $\zn$,
so the quadratic has two roots in~$\zn$.
The product of the roots is $3$, which is a non-square; so one root
(say~$x_1$) has odd order, and the other (say $x_2$) has even order.
We also note that the values $y_1=x_1+1$ and $y_2=x_2+1$ satisfy the 
quadratic equation $y^2+y+1=0$, and so $\ord_n(y_1)=\ord_n(y_2)=3$.

\begin{theorem}
Let $n$ be a prime satisfying $n \equiv 7$ or $31$ \mod{36},
$n>7$.   Suppose that the elements $x_1$ and $x_2$ from $\un$ that satisfy
$x^2 + 3x + 3 \equiv 0$ \mod{n} are such that
$\ord_n(x_1) = (n-1)/6$. Then $\ord_n(-(x_1 +2 )) = 3$ and so
\[ \un = \cy{-x_1 - 2}{3} \times \cy{-1}{2} \times \cy{x_1}{m} \]
where $m = (n-1)/6$.
\label{t21}
\end{theorem}

\begin{pf}
As noted above, $3 = \ord_n(x_2 + 1)$, and $x_2+1=-x_1-2$, since
$x_1+x_2=-3$.
\end{pf}

\paragraph{Coverage}
In the range $n < 1000$, Theorem~\ref{t21} covers values as follows:
\renewcommand{\arraystretch}{1.4}
\begin{eqnarray*}
n =  31: && \cy{25}{3} \times \cy{30}{2} \times \cy{4}{5}     \\
n =  43: && \cy{6}{3} \times \cy{42}{2}  \times \cy{35}{7}    \\
n =  79: && \cy{55}{3} \times \cy{78}{2} \times \cy{22}{13}   \\
n = 211: && \cy{196}{3} \times \cy{210}{2} \times \cy{13}{35} \\ 
n = 463: && \cy{21}{3} \times \cy{462}{2} \times \cy{440}{77} \\
n = 571: && \cy{109}{3} \times \cy{570}{2} \times \cy{460}{95}\\
n = 751: && \cy{678}{3} \times \cy{750}{2} \times \cy{71}{125}\\
n = 907: && \cy{522}{3} \times \cy{906}{2} \times \cy{383}{151}
\end{eqnarray*}

\begin{theorem}
Let $n$ be a prime satisfying $n \equiv 7$ or $31$ \mod{36},
$n>7$.   Suppose that the elements $x_1$ and $x_2$ from $\un$ that
satisfy $x^2 + 3x + 3 \equiv 0$ \mod{n} are such that
$\ord_n(x_1) = (n-1)/2$ and $\ord_n(x_2) = (n-1)$.   Then
$\ord_n(x_2 + 1) = 3$ and $\ord_n(2x_2 + 3) = (n-1)/6$, so that
\[ \un = \cy{2x_2 + 3}{m} \times \cy{x_2 + 1}{3} \times \cy{-1}{2} \]
where $m = (n-1)/6$.
\label{t22}
\end{theorem}

\begin{pf} 
Since $\un=\langle-1\rangle_2\times\langle x_1+1\rangle\times C_{(n-1)/6}$,
the hypothesis $\ord_n(x_2)=n-1$ shows that 
$x_2=-(x_1+1)c$ or $-(x_2+1)c$, where $\ord_n(c)=(n-1)/6$. 
Hence either $-(x_1+1)x_2$ or $=-(x_2+1)x_2$ has order $(n-1)/6$. 
Now $x_1+x_2=-3$ and $x_1x_2=3$, so $-(x_1+1)x_2=-3-x_2=x_1$, 
which has order $(n-1)/2$, by assumption. So 
\[-(x_2+1)x_2=-x_2^2-x_2=2x_2+3\]
has order $(n-1)/6$.
\end{pf}

\paragraph{Coverage}
In the range $n < 1000$, Theorem~2.2 covers values as follows:
\begin{eqnarray*}
n =  67: && \cy{59}{11} \times \cy{29}{3} \times \cy{66}{2}     \\
n = 103: && \cy{10}{17} \times \cy{46}{3} \times \cy{102}{2}    \\
n = 151: && \cy{86}{25} \times \cy{118}{3} \times \cy{150}{2}   \\
n = 367: && \cy{200}{61} \times \cy{283}{3} \times \cy{366}{2}  \\
n = 439: && \cy{343}{73} \times \cy{171}{3} \times \cy{438}{2}  \\
n = 499: && \cy{279}{83} \times \cy{139}{3} \times \cy{498}{2}  \\
n = 619: && \cy{505}{103} \times \cy{252}{3} \times \cy{618}{2} \\
n = 643: && \cy{355}{107} \times \cy{177}{3} \times \cy{642}{2} \\
n = 727: && \cy{563}{121} \times \cy{281}{3} \times \cy{726}{2} \\
n = 787: && \cy{28}{131} \times \cy{407}{3} \times \cy{786}{2}  \\
n = 967: && \cy{682}{162} \times \cy{824}{3} \times \cy{966}{2} \\
\end{eqnarray*}

\begin{theorem}
Let $n$ be a prime satisfying $n \equiv 7$ or $31$ \mod{36}.
Suppose that~$z$ is one of the elements $x_1$ and $x_2$ that satisfy
$x^2 + 3x + 3 \equiv 0$ \mod{n} and that
$\ord_n(2^{-1} z) = (n-1)/6$.   Then 
\[ \un = \cy{z+1}{3} \times \cy{2^{-1}z}{m} \times \cy{-1}{2} \]
where $m = (n-1)/6$.
\label{t23}
\end{theorem}

\begin{pf} As for Theorem~2.1.   But it depends on $n$
whether $z$ is the solution of $x^2 + 3x + 3 \equiv 0$ that has the 
larger or smaller order, and whether $z$ is $x_1$ or $x_2$.
\end{pf}

\paragraph{Coverage}
In the range $n < 1000$, Theorem~2.3 covers values as follows:
\begin{center}
$
\begin{array}{rcl}
n =  31: & \cy{5}{3} \times \cy{2}{5} \times \cy{30}{2}       & (z = x_1) \\
n =  67: & \cy{29}{3} \times \cy{14}{11} \times \cy{66}{2}    & (z = x_2) \\
n = 103: & \cy{56}{3} \times \cy{79}{17} \times \cy{102}{2}   & (z = x_1) \\
n = 151: & \cy{32}{3} \times \cy{91}{25} \times \cy{150}{2}   & (z = x_1) \\
n = 211: & \cy{196}{3} \times \cy{203}{35} \times \cy{210}{2} & (z = x_2) \\
n = 283: & \cy{44}{3} \times \cy{163}{47} \times \cy{282}{2}  & (z = x_2) \\
n = 691: & \cy{437}{3} \times \cy{218}{115} \times \cy{690}{2}& (z = x_2) \\
n = 787: & \cy{407}{3} \times \cy{203}{131} \times \cy{786}{2}& (z = x_2) \\
n = 823: & \cy{648}{3} \times \cy{735}{137} \times \cy{822}{2}& (z = x_1) \\
n = 907: & \cy{522}{3} \times \cy{714}{151} \times \cy{906}{2}& (z = x_2) \\
\end{array}
$
\end{center}

\paragraph{Note 2.1} In the range $n < 1000$,
Theorems~\ref{t21}--\ref{t23} exclude $n = 139$, 223, 331, 547, 607 and 859.
All but one of these has $x$-values $x_1$ and $x_2$ with
$\ord_n(x_1) = (n-1)/2$ and $\ord_n(x_2) = (n-1)/3$; the exception is
$n = 547$, which has $x_1 = 505$ and $x_2 = 39$, with
$\ord_n(x_1) = (n-1)/26$ and $\ord_n(x_2) = (n-1)/13$.

\paragraph{Problem 3} It is natural to wonder whether there are infinitely
many primes for which the conditions of one of the above theorems are 
satisfied. Here are some thoughts on this. In all cases we seek primes
congruent to $7$ or $31$ \mod{36}; Dirichlet's Theorem \cite[Theorem~2.10]{JJ}
guarantees that infinitely many such primes exist, and indeed they have 
density $1/6$ among all primes.

Consider Theorem~\ref{t21}. We require that an element of order $(n-1)/6$
(necess\-arily a sixth power) should satisfy $x^2+3x+3=0$, so there should 
be a solution $y$ of the equation $y^{12}+3y^6+3=0$. The Chebotarev density
theorem \cite{Chebotarev}, \cite[section 1.2.2]{Chebotarev2} 
guarantees that this equation will 
have a solution in an infinite set (indeed, a set of positive density) of 
primes. This theorem can further guarantee a set of positive density for 
which $x_1$ has six distinct sixth roots (so that $n\equiv1$ \mod{6}) and 
$x_2$ is a non-square (so that $n\equiv3$ \mod{4}), but we do not know how 
to exclude $n\equiv19$ \mod{36}.

A more serious difficulty is that the fact that $x$ is a sixth power 
guarantees only that its order divides $(n-1)/6$; it does not seem
easy to show that the order is precisely this value. Clearly this would 
be the case if $n=6q+1$ with $q$ prime; but it is not even known whether
infinitely many primes of this form occur.

Of the $1614$ primes less than $10^5$ which are congruent to $7$ or $31$
\mod{36}, there are $494$, $476$ and~$476$ that satisfy the conditions of 
Theorems~\ref{t21}--\ref{t23} respectively 

\medbreak

We conclude this subsection with a converse to the preceding theorems.

\begin{theorem}
Any \textup{3AP} decomposition of $\un$ for $n$ prime, in which the 
generators have orders $2$, $3$ and $(n-1)/6$, arises as in one of 
the three preceding theorems.
\label{t24}
\end{theorem}

\begin{pf}
We already saw that $n$ must be congruent to $7$ or $31$ \mod{36}. The
only element of order~$2$ is $-1$, and the only elements of order~$3$ are
$x_1+1$ and $x_2+1$, where $x_1$ and $x_2$ are the roots of $x^2+3x+3=0$
(with the convention that $x_1$ has odd order and $x_2$ even order).
The only possibilities for the third generator are thus
$-x_i-3$, $2x_i+3$, or $2^{-1}x_i$, for $i=1$ or $i=2$. We treat the three
cases in turn.

In the first case, since $x_1+x_2=-3$, we have $-x_1-3=x_2$ and 
\emph{vice versa}. Since $x_2$ has even order by our convention, we must
have $i=2$, and the generators are $x_1$, $-1$ and $x_2+1$; the requirement
is that $x_1$ has order $(n-1)/6$.

In the second case, we assume that $2x_i+3$ has order $(n-1)/6$, and have to
prove that $i=2$ and that the orders of $x_1$ and $x_2$ are $(n-1)/2$ and
$(n-1)$ respectively. Let $j=3-i$. From the proof of Theorem~\ref{t22}, we see
that $2x_i+3=-(x_i+1)x_i$, so that $x_i=-(x_j+1)(2x_i+3)$, the product of
elements of orders $2$, $3$ and $(n-1)/6$; so $x_i$ has order $n-1$. Thus
$i=2$. Now $(x_1+1)x_1=2x_2+3$ has order $(n-1)/6$, so $x_1=(x_2+1)(2x_2+3)$
has order $(n-1)/2$.

Finally, the third case obviously gives the situation of Theorem~\ref{t23}.
\end{pf}

\subsection{The case $n-1 = 3 \cdot 4 \cdot \mu$}
\label{ss22}

We now prove two theorems that apply for prime values $n$ such
that the factors in a 3AP decomposition of $\un$ have orders
$3$, $4$ and $\mu$ where $3$, $4$ and $\mu$ are pairwise co-prime.
We give no theorem for the situation where $\mu$ is the order of the
middle generator. This case can occur; the smallest example is for $n=997$.

\begin{theorem}
Let $n$ be a prime satisfying $n \equiv 13$, $61$, $85$ or $133$
\mod{144}, $n > 13$.   Suppose that there is an element $x$ from
$\un$ such that $x^2 + 3x + 3 \equiv 0$ \mod{n} and such
that there is also an element $k$ with $\ord_n(x+1+k) = 4$ and
$\ord_n(x+1+2k) = (n-1)/12$.   Then
\[ \un = \cy{x+1}{3} \times \cy{x+1+k}{4} \times \cy{x+1 +2k}{\mu}  \]
where $\mu = (n-1)/12$.
\label{t25}
\end{theorem}

\begin{pf} 
As for Theorem~2.1. Note that the condition on $x+1+k$
can be written $(x+1+k)^2\equiv-1$ \mod{n}.
\end{pf}

\paragraph{Coverage}
In the range $n < 1000$, Theorem~\ref{t25} covers values as follows:
\begin{center}
$
\begin{array}{rcl}
 n & \textup{3AP decomposition of }\un & \ord_n(x) \\
\hline
 61 & \langle  13 \rangle_3 \times \langle  11 \rangle_4
                         \times \langle   9 \rangle_5      & 15 = (n-1)/4  \\
349 & \langle 122 \rangle_3 \times \langle 213 \rangle_4
                         \times \langle 304 \rangle_{29}   & 58 = (n-1)/6  \\
661 & \left\{ \begin{array}{@{}l@{}}
      \langle 364 \rangle_3 \times \langle 106 \rangle_4
                         \times \langle 509 \rangle_{55}       \\
      \langle 364 \rangle_3 \times \langle 555 \rangle_4
                         \times \langle  85 \rangle_{55}
      \end{array} \right\}                                 & 66 = (n-1)/10 \\
\hline
\end{array}
$
\end{center} 
For the two examples for $n = 661$, the generators of order 55 are related
by the congruence $85 \equiv 509^3$ \mod{661}.

\begin{theorem}
Let $n$ be a prime satisfying $n \equiv 13$, $61$, $85$ or $133$
\mod{144}, $n > 13$.   Suppose that there is an element $x$ from
$\un$ such that $x^2 + 3x + 3 \equiv 0$ \mod{n} and such
that there is also an element $k$ with $\ord_n(x+1+k) = 4$ and
$\ord_n(x+1-k) = (n-1)/12$.   Then
\[ \un = \langle x+1-k \rangle_\mu \times \langle x+1 \rangle_3 
                                   \times \langle x+1+k \rangle_4         \]
where $\mu = (n-1)/12$.
\label{t26}
\end{theorem}

\begin{pf} As for Theorem~2.1.
\end{pf}

\paragraph{Coverage}
In the range $n < 1000$, Theorem~\ref{t26} covers values as follows:
\begin{center}
$
\begin{array}{rcl}
 n & \textup{3AP decomposition of }\un & \ord_n(x) \\
\hline
157 & \langle 153 \rangle_{13} \times \langle  12 \rangle_3
                               \times \langle  28 \rangle_4 &  39 = (n-1)/4 \\
229 & \langle 161 \rangle_{19} \times \langle 134 \rangle_3
                               \times \langle 107 \rangle_4 & 228 = (n-1)   \\
349 & \langle  31 \rangle_{29} \times \langle 122 \rangle_3
                               \times \langle 213 \rangle_4 &  58 = (n-1)/6 \\
373 & \langle  91 \rangle_{31} \times \langle 284 \rangle_3
                               \times \langle 104 \rangle_4 &  93 = (n-1)/4 \\
997 & \langle 226 \rangle_{83} \times \langle 692 \rangle_3
                               \times \langle 161 \rangle_4 & 498 = (n-1)/2 \\
\hline
\end{array}
$
\end{center}

\medskip

In the range $n < 1000$, Theorems \ref{t25} and \ref{t26} fail to provide
3AP decompositions of $\un$ for $n = 277$, 421, 709, 733, 853 and 877.
However, 3AP decompositions for $n = 421$ exist for other partitions of $n-1$.

\subsection{The case $n-1 = 2 \cdot 5 \cdot \nu$}
\label{ss23}

We now examine what occurs for primes $n$ such that the factors in a
3AP decomposition of $\un$ have orders 2, 5 and $\nu$, these orders
being pairwise co-prime.
Amongst primes satisfying $n \equiv 11$, 31, 71 and 91 (mod~100), 
$n > 11$, the patterns of occurrence of such 3AP decompositions are 
very similar to those reported in~\S\ref{ss21} above. 
For each relevant value of $n$ there are 4 elements of order 5; their sum
is ${-1}$ and their product is ${+1}$.   Sometimes more than one of
the four can be used.

\paragraph{Type 2.3(a)}
Analogous to the decompositions obtainable via Theorem \ref{t21}, we now
have 3AP decompositions of the form 
\[ \un = \cy{-z-2}{5} \times \cy{-1}{2} \times \cy{z}{\nu}\ .          \]
In the range $n < 1000$ they are as follows:
\begin{eqnarray*}
 n =  31: && \cy{4}{5} \times \cy{30}{2} \times \cy{25}{3}             \\
 n = 191: &  \left\{ \begin{array}{@{}c@{}} \\ \\ \end{array} \right.
           & \begin{array}{@{}l@{}}
             \cy{39}{5} \times \cy{190}{2} \times \cy{150}{19} \\
             \cy{184}{5} \times \cy{190}{2} \times \cy{5}{19}  \\
             \end{array}                                               \\
 n = 271: &  \left\{ \begin{array}{@{}c@{}} \\ \\ \end{array} \right.
           & \begin{array}{@{}l@{}}
             \cy{10}{5} \times \cy{270}{2} \times \cy{259}{27} \\
             \cy{244}{5} \times \cy{270}{2} \times \cy{25}{27} \\
             \end{array}                                               \\
 n = 431: && \cy{405}{5} \times \cy{430}{2} \times \cy{24}{43}         \\
 n = 691: && \cy{89}{5} \times \cy{690}{2} \times \cy{600}{69}         \\
 n = 991: && \cy{799}{5} \times \cy{990}{2} \times \cy{190}{99}
\end{eqnarray*}

\paragraph{Type 2.3(b)}
Analogous to the decompositions obtainable via Theorem \ref{t22}, we 
have 3AP decompositions of the form 
\[ \un = \cy{2z+1}{\nu} \times \cy{z}{5} \times \cy{-1}{2}\ .          \]
In the range $n < 1000$ they are as follows:
\begin{eqnarray*}
 n =  31: && \cy{5}{3} \times \cy{2}{5} \times \cy{30}{2}              \\
 n = 131: && \cy{107}{13} \times \cy{53}{5} \times \cy{130}{2}         \\
 n = 311: &  \left\{ \begin{array}{@{}c@{}} \\ \\ \end{array} \right.
           & \begin{array}{@{}l@{}}
             \cy{13}{31} \times \cy{6}{5} \times \cy{310}{2}   \\
             \cy{105}{31} \times \cy{52}{5} \times \cy{310}{2} \\
             \end{array}                                               \\
 n = 491: && \cy{203}{49} \times \cy{101}{5} \times \cy{490}{2}        \\
 n = 811: && \cy{330}{81} \times \cy{570}{5} \times \cy{810}{2}        \\
 n = 991: && \cy{395}{99} \times \cy{197}{5} \times \cy{990}{2}
\end{eqnarray*} 

\paragraph{Type 2.3(c)}
Analogous to the decompositions obtainable via Theorem \ref{t23}, we 
have 3AP decompositions of the form 
\[ \un = \cy{2z+1}{5} \times \cy{z}{\nu} \times \cy{-1}{2}\ .          \] 
In the range $n < 1000$ they are as follows:
\begin{eqnarray*}
 n = 271: && \cy{10}{5} \times \cy{140}{27} \times \cy{270}{2}         \\
 n = 691: && \cy{132}{5} \times \cy{411}{69} \times \cy{690}{2}        \\
 n = 971: && \cy{803}{5} \times \cy{401}{97} \times \cy{970}{2}        \\
 n = 991: && \cy{197}{5} \times \cy{98}{99} \times \cy{990}{2}
\end{eqnarray*}

\paragraph{Note 2.2} In the range $n < 1000$, no 3AP decomposition of any
of the types 2.3(a), 2.3(b) or 2.3(c) exists for $n = 71$, $211$, $331$, 
$571$, $631$ or $911$.

\subsection{Some double-barrelled cases}
\label{ss24}

Two special cases arise for prime $n$ such that $\un$ has more than
one 3AP decomposition.   These are where we can write either
\[ \un = \cy{k}{a} \times \cy{k+z}{b} \times \cy{k+2z}{c}
       = \cy{k-2z}{b} \times \cy{k}{a} \times \cy{k+2z}{c}         \]
or
\[ \un = \cy{k}{a} \times \cy{k+z}{b} \times \cy{k+2z}{c}
       = \cy{k+z}{b} \times \cy{k+2z}{c} \times \cy{k+3z}{a}\ .    \]
For the first of these we need $k+z$ and $k - 2z$ to have the same
order \mod{n} and each to be a power of the other \mod{n}.   
For the second we need the same relationship between $k$ and $k+3z$.   
We have failed to find any theorems to indicate when these cases arise.
In the range $n < 1000$, the occurrences of the first case are these:
\begin{center}
$
\begin{array}{lcccc}
 \uu{67}  & = & \cy{29}{3} \times \cy{14}{11} \times \cy{66}{2}
          & = & \cy{59}{11} \times \cy{29}{3} \times \cy{66}{2}         \\
 \uu{211} & = & \cy{210}{2} \times \cy{203}{35} \times \cy{196}{3}
          & = & \cy{13}{35} \times \cy{210}{2} \times \cy{196}{3}       \\
 \uu{271} & = & \cy{270}{2} \times \cy{140}{27} \times \cy{10}{5}
          & = & \cy{259}{27} \times \cy{270}{2} \times \cy{10}{5}       \\
 \uu{331} & = & \cy{167}{11} \times \cy{83}{15} \times \cy{330}{2}
          & = & \cy{4}{15} \times \cy{167}{11} \times \cy{330}{2}       \\
 \uu{379} & = & \cy{378}{2} \times \cy{119}{7} \times \cy{239}{27}
          & = & \cy{138}{7} \times \cy{378}{2} \times \cy{239}{27}      \\
 \uu{661} & = & \cy{364}{3} \times \cy{391}{20} \times \cy{418}{11}
          & = & \cy{310}{20} \times \cy{364}{3} \times \cy{418}{11}     \\
 \uu{787} & = & \cy{407}{3} \times \cy{203}{131} \times \cy{786}{2}
          & = & \cy{28}{131} \times \cy{407}{3} \times \cy{786}{2}      \\
 \uu{907} & = & \cy{906}{2} \times \cy{714}{151} \times \cy{522}{3}
          & = & \cy{383}{151} \times \cy{906}{2} \times \cy{522}{3}     \\
\end{array}
$
\end{center}
whereas the occurrences of the second case are these:
\begin{center}
$
\begin{array}{lcccc}
 \uu{349} & = & \cy{31}{29} \times \cy{122}{3} \times \cy{213}{4}
          & = & \cy{122}{3} \times \cy{213}{4} \times \cy{304}{29}      \\
 \uu{599} & = & \cy{578}{23} \times \cy{598}{2} \times \cy{19}{13}
          & = & \cy{598}{2} \times \cy{19}{13} \times \cy{39}{23}\ .    \\
\end{array}
$
\end{center}

\section {Lifts}
\label{s3}

We now consider how and when an AP decomposition of $\uu{n}$ can be used 
to obtain AP decompositions of $\uu{n'}$ where $n'$ is a power of $n$ or
some other multiple of~$n$.

\subsection{Definitions}
\label{ss31}

In this section, we allow \emph{weak} 3AP decompositions 
$\uu{n}=\cy{x}{a}\times\cy{y}{b}\times\cy{z}{c}$,
where one of $x,y,z$ is allowed to be $1$ (so that the corresponding 
cyclic factor is trivial). For example,
\begin{eqnarray}
\uu{7}=\cy{4}{3}\times\cy{6}{2}\times\cy{1}{1}
\label{sct31}
\end{eqnarray}
and
\[\uu{103}=\cy{1}{1}\times\cy{47}{6}\times\cy{93}{17}\]
are weak 3AP decompositions. Where it aids clarity, we refer to a 3AP
decomposition in the original sense as being \emph{strong}.

If $p$ is a prime satisfying $p \equiv 11$ \mod{12} and $\ord_p(3) = (p-1)/2$,
then $\uu{p} = \cy{-1}{2}\times\cy{1}{1}\times\cy{3}{(p-1)/2}$.
Likewise if $p$ is a prime satisfying $p \equiv 7$ \mod{12} and
$\ord_p(-3) = (p-1)/2$, then
$\uu{p} = \cy{-3}{(p-1)/2}\times\cy{-1}{2}\times\cy{1}{1}$.

If $n$ divides $n'$, then the map $x\mapsto x$ \mod{n} is a ring epimorphism
from $\zz{n'}$ to $\zz{n}$, and maps $\uu{n'}$ onto $\uu{n}$. It preserves
the property of forming an arithmetic progression. However, it does not
in general map a (weak) 3-AP decomposition of $\uu{n'}$ to a (weak)
3-AP decomposition of $\uu{n}$. (It maps a generating set to a generating
set, but does not necessarily preserve the direct sum decomposition.)
Note in passing that, if $n$ divides $n'$, then $\phi(n)$ divides $\phi(n')$,
the quotient being the order of the kernel of the homomorphism from
$\uu{n'}$ to $\uu{n}$.

Suppose that $n$ divides $n'$, and that the 3AP decompositions
\[\uu{n}=\cy{x}{a}\times\cy{y}{b}\times\cy{z}{c}\]
and 
\[\uu{n'}=\cy{x'}{a'}\times\cy{y'}{b'}\times\cy{z'}{c'}\]
satisfy 
$x \equiv x'$ \mod{n}, $y' \equiv y$ \mod{n} and $z \equiv z'$ \mod{n}. 
Then we call the
second decomposition a \emph{lift} of the first, with \emph{index} $n'/n$. 
Note that we must have $a\mid a'$, $b\mid b'$, $c\mid c'$, and 
$(a'b'c')/(abc)=\phi(n')/\phi(n)$. For example, the decomposition 
\begin{eqnarray}
\uu{49}=\cy{18}{3}\times\cy{48}{2}\times\cy{29}{7}
\label{sct32}
\end{eqnarray}
is a lift of the weak 3AP decomposition (\ref{sct31}) of $\uu{7}$.
We further describe $x'$ as being a \emph{lift} of $x$, and so on.

\subsection{Lifts from $n$ to $np$, with $p$ an odd prime}
\label{ss32}

We are unable to give necessary and sufficient conditions for lifts to
exist. In the remainder of this section, we consider the case where $n'=np$
for an odd prime $p$. (We do not know whether every lift can be obtained
by a sequence of lifts where the indices are primes.)

We subdivide the analysis into three cases.

\paragraph{Case 1: $p^2$ divides $n$.} Let $n=p^km$ ($k \ge 2$) where 
$p$ does not divide $m$.

In this case, if an element $x\in\uu{n}$ has $\ord_{p^k}(x)$ divisible by $p$,
then any lift $x'$ of $x$ satisfies $\ord_{n'}(x')=p\,\ord_n(x)$. So at most
one of $x,y,z$ can satisfy this condition if there is a lift which is a
3AP decomposition. Conversely, if, say, $x$~satisfies the condition but
$y$ and $z$ do not, then we can choose lifts $y'$
and $z'$ of $y$ and $z$ satisfying $\ord_{n'}(y')=\ord_n(y)$ and
$\ord_{n'}(z')=\ord_n(z)$, and the unique lift $x'$ of $x$ such that
$(x',y',z')$ is an AP in $\zz{n'}$, to obtain a 3AP decomposition 
of $\uu{n'}$.

For example, below (Case 2, Subcase 2.2) we find a decomposition
\[\uu{275}=\cy{274}{2} \times \cy{166}{5} \times \cy{58}{20}.\]
This cannot be lifted to a 3AP decomposition of $\uu{n'}$ with
$n' = 5 \times 275$. On the other hand, the decomposition (\ref{sct32})
can be lifted to
\[\uu{343}=\cy{18}{3}\times\cy{342}{2}\times\cy{323}{49}.\]

\paragraph{Case 2:} $p$ divides $n$ but $p^2$ does not. We further subdivide
into four cases according to how many of the orders $a,b,c$ are divisible
by~$p$. Note that $\phi(np)=\phi(n)p$.

For any
$x\in\uu{n}$, with $\ord_n(x)=a$, the lifts of $x$ to $\uu{pn}$ belong to
an extension of $C_p$ by $C_a$, which is isomorphic to $C_p\times C_a$. 
Hence, if $p$ does not divide $a$, then one of these lifts (which
we call the \emph{special lift}) has order~$a$, and the other $p-1$ have
order~$pa$; while if $p$ divides $a$, then all have order~$a$.

\subparagraph{Subcase 2.1:} None of $a,b,c$ is divisible by~$p$. If a lift
is a 3AP decomposition, then two of $x,y,z$ lift to elements of the same
orders (and so must be special lifts), while the third lifts to an element
with $p$ times the order. Let
$x',y',z'$ be the special lifts of $x,y,z$. We call the decomposition 
\emph{unproductive} if $(x',y',z')$ is an AP in $\zz{np}$. In this case,
there is no lift to a 3AP decomposition of $\uu{np}$. In the contrary 
\emph{productive} case, there are three lifts which are 3AP decompositions 
since we may choose for which two of $x,y,z$ we use special lifts, and the 
third lift is determined by the AP requirement.

Thus, if we start from the productive strong 3AP decomposition
\[\uu{31} = \cy{25}{3} \times \cy{30}{2} \times \cy{4}{5}\]
we obtain the three lifts
\begin{eqnarray*}
\uu{31^2} 
&=& \cy{521}{3} \times \cy{960}{2} \times \cy{438}{155} \\
&=& \cy{521}{3} \times \cy{526}{62} \times \cy{531}{5} \\
&=& \cy{428}{93} \times \cy{960}{2} \times \cy{531}{5}\,.
\end{eqnarray*}
Likewise we can start from the productive strong 3AP decomposition
\begin{eqnarray*}
 \uu{35}  & = & \cy{11}{3} \times \cy{34}{2} \times \cy{22}{4}\ .
\end{eqnarray*}
to obtain the lifts
\begin{eqnarray*}
 \uu{245} & = & \cy{116}{3} \times \cy{244}{2} \times \cy{127}{28}   \\
          & = & \cy{116}{3} \times \cy{34}{14} \times \cy{197}{4}    \\
          & = & \cy{46}{21} \times \cy{244}{2} \times \cy{197}{4}
\end{eqnarray*}
and
\begin{eqnarray*}
 \uu{175} & = & \cy{151}{3} \times \cy{174}{2} \times \cy{22}{20}    \\
          & = & \cy{151}{3} \times \cy{104}{10} \times \cy{57}{4}    \\
          & = & \cy{116}{15} \times \cy{174}{2} \times \cy{57}{4}\ .
\end{eqnarray*}

If we start from a productive weak 3AP decomposition, then two of 
the three lifts are weak but the third (where the identity lifts to an 
element of order~$p$) is strong. This happens in the example (\ref{sct32})
of a 3AP decomposition of $\uu{49}$; the two corresponding weak 3AP 
decompositions lifted from (\ref{sct31}) are
\[\uu{49}=\cy{18}{3}\times\cy{34}{14}\times\cy{1}{1}
=\cy{46}{21}\times\cy{48}{2}\times\cy{1}{1}.\]

Now consider the prime $n = 379$.  The strong 3AP decomposition
\begin{eqnarray}
\uu{379} &=& \cy{239}{27}\times\cy{378}{2}\times\cy{138}{7}
\label{unpro1}
\end{eqnarray}
is unproductive: the special lifts of $239$, $378$ and $138$ are
respectively $8956$, $143640$ and $134683$, which happen to be in
arithmetic progression \mod{379^2}. The two weak 3AP decompositions
\begin{eqnarray}
\uu{11} &=& \cy{10}{2}\times\cy{1}{1}\times\cy{3}{5}
\label{unpro2}
\end{eqnarray}
and
\begin{eqnarray}
\uu{461} &=& \cy{1}{1}\times\cy{48}{4}\times\cy{95}{115}
\label{unpro3}
\end{eqnarray}
are also unproductive, but these are the only ones with prime modulus
less than $1000$. It appears that productive decompositions predominate;
unproductive ones depend on an accidental coincidence which is comparatively
rare.

\subparagraph{Subcase 2.2:} One of $a,b,c$ (say $a$, without loss) is
divisible by $p$. Now choose the special lift of either $y$ or $z$, and any
non-special lift of the other; the lift of $x$ is determined by the AP
requirement. So there are $2(p-1)$ lifts to 3AP decompositions.

Here is an example. Start from a weak 3AP decompositions of $\uu{55}$:
\begin{eqnarray}
\uu{55}=\cy{54}{2}\times\cy{1}{1}\times\cy{3}{20}\,.
\label{sct33}
\end{eqnarray}
We wish to lift to strong 3AP decompositions of $\uu{275}$. We are 
in this subcase. All lifts of $3$ have order~$20$, but each of the
generators $54$ and $1$ has one special lift (namely $274$ and $1$ 
respectively). So we must use a non-special lift of $1$, the special 
lift of $54$, and the lift of $3$ which completes the AP:
\begin{eqnarray*}
 \uu{275} & = & \cy{274}{2} \times \cy{166}{5} \times \cy{58}{20}   \\
          & = & \cy{274}{2} \times \cy{56}{5} \times \cy{113}{20}   \\
          & = & \cy{274}{2} \times \cy{221}{5} \times \cy{168}{20}  \\
          & = & \cy{274}{2} \times \cy{111}{5} \times \cy{223}{20}
\end{eqnarray*}
The other four lifts (where we use the special lift of $1$ and a non-special
lift of~$54$) are weak 3AP decompositions. In the same way, the decomposition
\begin{eqnarray}
\uu{55} = \cy{52}{20} \times \cy{54}{2} \times \cy{1}{1}
\label{sct34}
\end{eqnarray}
gives rise to four more strong 3AP decompositions of $\uu{275}$.

Suppose, however, that we consider lifting (\ref{sct33}) and (\ref{sct34})
from $n = 55$ to \mbox{$n = 605$}.   The only lifts of 52, 54, 1 and 3
\mod{605} which have orders 20, 2, 1, and 20 respectively are
602, 604, 1 and 3 \mod{605}, which are in AP, so we fail to obtain
any 3AP decomposition for $\uu{605}$. 

\subparagraph{Subcase 2.3:} Two of $a,b,c$ (say $a$ and $b$) are divisible
by $p$. We must choose a non-special lift of $z$, and any lift of $x$; so
there are $p(p-1)$ lifts to 3AP decompositions of~$\uu{n'}$.   Suppose,
for example, that we take $n = 273 = 3\times7\times13$ and $p=3$, to give
$n' = 819$. Computer enumeration has shown that there are $108$ strong
3AP decompositions of~$\uu{273}$, each perforce having 2 generators whose
orders are multiples of~3.   The 648 lifts to $\uu{819}$ arise from the
strong decompositions.

\subparagraph{Subcase 2.4:} All three of $a,b,c$ are divisible by $p$. In
this case, three distinct prime divisors of $n$ are congruent to $1$ \mod{p}, 
and so $np$ has at least four prime divisors, so no
3AP decomposition of $\uu{np}$ can exist. (Alternatively, note that all
lifts of $x,y,z$ have the same orders as the original elements, so the group
they generate has the same order as $\langle x,y,z\rangle$.)

\paragraph{Case 3:} $p$ does not divide~$n$. In this case, 
$\uu{np}\cong\uu{n}\times\uu{p}$. This case is the most difficult and 
we do not have any general criteria for a lift to exist. However, 
$\uu{np}$ is a product of at most 3 cyclic groups.   From the structure
of the group of units, as described in the Introduction, we see that
$n=q^\alpha r^\beta$, or $2q^\alpha r^\beta$,
or $4q^\alpha$, or $2^\alpha$, for some odd primes $q$ and $r$, and some
$\alpha,\beta\ge0$. In this case, one of the lifts of any generator of a
3AP decomposition of $\uu{n}$ is a multiple of $p$, and therefore must be
disallowed as a \emph{spurious lift}. The spurious lifts $a',b',c'$ of the 
three generators $a,b,c$ are in~AP. For they form an AP (mod~$n$) by
definition, and a trivial AP (mod~$p$). By the Chinese Remainder Theorem,
the congruences
\[b'-a'\equiv c'-b'\hbox{ (mod~$n$)},\qquad b'-a'\equiv c'-b'\hbox{ (mod~$p$)}\]
imply that $b'-a'\equiv c'-b'$ (mod~$pn$).

Suppose that we have $n = 31$ and $p = 5$, and we consider lifting
\[\uu{31} = \cy{25}{3}\times\cy{30}{2}\times\cy{4}{5}\,.\]
The respective spurious lifts are 25, 30 and 35.   Two lifts of the
3AP decomposition are available:
\begin{eqnarray*}
\uu{155} & = & \cy{56}{3}\times\cy{154}{2}\times\cy{97}{20}\\
         &   & \cy{87}{12}\times\cy{154}{2}\times\cy{66}{5}\,.
\end{eqnarray*}

\subsection{Lifting to $\uu{p^\alpha}$}
\label{ss33}

In one special case of lifts we can draw a strong conclusion.
This case requires a productive 3AP decomposition for $\uu{p}$ where
$p$ is prime; we recall the unproductive examples (\ref{unpro1}),
(\ref{unpro2}) and (\ref{unpro3}) given above. 

\begin{theorem}
Let $p$ be an odd prime, and suppose that
\[\uu{p}=\cy{x}{a}\times\cy{y}{b}\times\cy{z}{c}\]
is a productive (possibly weak) \textup{3AP} decomposition. Then for any
$\alpha\ge2$, there is a lift of the given decomposition which is a (strong)
\textup{3AP} decomposition of~$\uu{p^\alpha}$.
\end{theorem}

\begin{pf}
We show inductively that there is a lift where two of the lifted elements
have the same orders as the originals, and the third has order multiplied
by~$p^{\alpha-1}$.

For the first step, we are in case 2, subcase 2.1; in this case we saw
that any productive decomposition has three lifts, at least one of which
is strong.

For the general step, we are in case 1, and we start with a decomposition
in which two of the elements have orders coprime to $p$; so the necessary
condition for this case is satisfied, and the lift exists.
\end{pf}

\section{$n$ an odd composite integer}
\label{s4}

\subsection{$n$ a multiple of~$3$}
\label{ss41}

\begin{theorem} 
A \textup{3AP} decomposition of $\uu{3p}$ does not exist for any prime $p$
with $p > 3$.
\label{t3p}
\end{theorem}

\begin{pf} Suppose that 
\[\uu{3p}=\langle a\rangle\times\langle a+d\rangle\times\langle a+2d\rangle.\]
Then $d$ is divisible by $3$, since otherwise one of $a$, $a+d$, $a+2d$ would
be a multiple of~$3$. 

If $a\equiv1$~(mod~$3$) then all three generators are congruent
to $1$~(mod~$3$), and so is every element in the group they generate, which
is not possible. On the other hand, if $a\equiv2$~(mod~$3$), then each of
the generators has even order (since it has even order in $\uu{3}$), and so
$C_2\times C_2\times C_2\le\uu{3p}=C_2\times C_{p-1}$, a contradiction.
\end{pf}

The next result has a similar proof; it is rather special but rules out one
particular type of 3AP decomposition.

\begin{theorem}
There is no \textup{3AP} decomposition of $\uu{3m}$ of the form
\[\uu{3m}=\langle a\rangle\times\langle a+m\rangle\times\langle a+2m\rangle.\]
\end{theorem}

\begin{pf}
The argument of the preceding theorem shows that $m$ is divisible by~$3$.
Since all the generators are congruent mod~$m$, 
and projection from $\uu{3m}$ to $\uu{m}$ is onto, we see that $\uu{m}$ must
be cyclic (and $a$ is a primitive root of $m$), so $m$ is of the form $p^t$,
or $2p^t$ (for some odd prime $p$), or $m=4$. So necessarily $p=3$. But
now the order of $a$ (mod~$m$) is $2\cdot3^{t-1}$, and the same goes for the
other generators as well. Their orders (mod~$3m$) are at least as large, so
we must have $(2\cdot3^{t-1})^3\le 2\cdot3^t$, which is impossible.
\end{pf}

\subsection{Products of three primes}
\label{ss42}

Theorem~\ref{t3p} does not rule out 3AP decompositions of $\uu{3pq}$, where
$p$ and $q$ are distinct primes, and indeed these do exist. In this case, 
a new phenomenon occurs: we can obtain new solutions from old. This works
more generally for the case where $n$ is the product of three odd primes
$p,q,r$, and $\xi(n)=4$ (so that $\uu{n}=C_{\lambda_n}\times C_2\times C_2$).

Suppose that the abelian group $A$ can be written (adapting our previous
notation) as
\[A=\cy{x}{2a}\times\cy{y}{2}\times\cy{z}{2}.\]
Then $A$ contains an elementary abelian group $B$ of of order~$8$ generated by
$x^a,y,z$. If three elements $x',y',z'$ have the properties that their orders
are $2a,2,2$ respectively and $\langle(x')^a,y',z'\rangle=B$, then $x',y',z'$
generate cyclic subgroups whose direct product is $A$.
If $A=\uu{n}$ for some $n$, then multiplying an
arithmetic progression by a fixed unit yields an arithmetic progression; so
we look for an element $u$ such that $x'=xu$, $y'=yu$, $z'=zu$ satisfy the
above conditions.

We see that $u$ must have order~$2$, so $u\in B$. If $a$ is even then 
$(xu)^a=x^a$, while if $a$ is odd then $(xu)^a=x^au$. It is then easy to
check that the allowable values of $u$ are as follows:
\begin{itemize}
\item $u\in\{x^a,yz,x^ayz\}$ if $a$ is even;
\item $u\in\{x^ay,x^az,yz\}$ if $a$ is odd.
\end{itemize}
In each case, the possible values of $u$, together with the identity, form
a subgroup of $B$; so no further expressions can be obtained by repeating
the procedure. Moreover, in each case, $yz$ is an allowed multiplier, and
converts $[x,y,z]$ into $[xyz,z,y]$; so the solutions come in pairs, each
pair consisting of the first three and the last three terms in the sequence
$[x,y,z,xyz]$.

\begin{theorem}
Suppose that $n$ is the product of three odd primes, and that
\[\uu{n}=\cy{x}{\lambda}\times\cy{y}{2}\times\cy{z}{2}\]
is a \textup{3AP} decomposition, where $\lambda=\lambda_n$. Then
\[\uu{n}=\cy{ux}{\lambda}\times\cy{uy}{2}\times\cy{uz}{2}\]
is also a \textup{3AP} decomposition, where
\begin{itemize} 
\item $u\in\{x^{\lambda/2},yz,x^{\lambda/2}yz\}$ if 
$\lambda\equiv0\,(\mathrm{mod}\,4)$;
\item $u\in\{x^{\lambda/2}y,x^{\lambda/2}z,yz\}$ if
$\lambda\equiv2\,(\mathrm{mod}\,4)$.
\end{itemize}
In each case, there are two four-term arithmetic progressions whose
three-term subprogressions give the stated decompositions.
\end{theorem}

We call these sets of four decompositions \emph{quartets}.
Here are some examples of quartets, in cases where one of the primes 
dividing $n$ is~$3$. We list the values of $n$ and $\lambda$, and the
two four-term progressions $[x,y,z,xyz]$; the orders of the terms are
$\lambda,2,2,\lambda$, and the first and last three give 3AP decompositions.

\paragraph{Case} $\lambda\equiv0$ (mod~$4$):
\begin{itemize}
\item $105$; $12$; $[38, 71, 104, 32]$, $[17, 29, 41, 53]$
\item $165$; $20$; $[113, 56, 164, 107]$, $[47, 89, 131, 8]$
\item $285$; $36$; $[98, 191, 284, 92]$, $[212, 134, 56, 263]$
\item $357$; $48$; $[122, 239, 356, 116]$, $[269, 50, 188, 326]$
\item $465$; $60$; $[158, 311, 464, 152]$, $[437, 404, 371, 338]$
\end{itemize}

\paragraph{Case} $\lambda\equiv2$ (mod~$4$):
\begin{itemize}
\item $231$; $30$; $[80, 155, 230, 74]$, $[179, 188, 197, 206]$
\item $483$; $66$; $[164, 323, 482, 158]$, $[95, 461, 344, 227]$
\end{itemize}

In general, there is no rquirement that an end-term in a quartet
should be the product \mod{n} of the other three terms.   A
counter-example is the following, where the subscript integers
are the orders of the terms:
\begin{itemize}
\item $315$; $12$; $[8_4, 131_6, 254_6, 62_4]$
\end{itemize}

%

\subsection{Some results for $n = pq$ $(p > 3, q > 3)$}

We now indicate how the role of the value ${-3}$ \mod{n}, as discussed
in \S\ref{s2} above, carries over to composite values of $n$.

\begin{theorem}
Let $p$ and $q$ be primes greater than $3$, with $p\equiv3$ \mod{4}, and 
suppose that $\ord_p(-3)=(p-1)/2$ and $\ord_q(-3)=q-1$. 
(This implies that $p\equiv1$ \mod{3} and $q\equiv2$ \mod{3}). 
Let $n=pq$, and let $x$ be the unique element of $\un$ congruent 
to $1$ \mod{p} and to $-3$ \mod{q}.  Then
\[\un = \langle-x-2\rangle_{(p-1)/2} \times \langle-1\rangle_2 \times
\langle x\rangle_{q-1}\,,\]
which is a lift of $\uu{p}=\cy{-3}{(p-1)/2}\times\cy{-1}{2}\times\cy{1}{1}$.
\label{t41}
\end{theorem}

\begin{pf}
The congruences \mod{3} arise by noticing that $-3$ is a quadratic residue
\mod{p} and non-residue \mod{q}, and applying quadratic reciprocity.

We have $\ord_p(x)=1$ and $\ord_q(x)=q-1$, so $\ord_n(x)=q-1$. Also, $-x-2$
is congruent to $-3$ \mod{p} and to $1$ \mod{q}, so $\ord_p(-x-2)=(p-1)/2$
and $\ord_q(-x-2)=1$, whence $\ord_n(-x-2)=(p-1)/2$.

Since $p\equiv3$ \mod{4}, we have $\uu{p} = \langle-3\rangle_{(p-1)/2}
\times \langle-1\rangle_2$, as the orders of the factors are co-prime.
So the group $A=\langle-x-2\rangle_{(p-1)/2} \times \langle-1\rangle_2 \times
\langle x\rangle_{q-1}$ projects onto $\mathbb{U}_p$. Also, $x$ belongs to
the kernel of this projection; since $x$ is a primitive root of $q$, the
kernel is $\mathbb{U}_q$. So $A=\un$.
\end{pf}

\paragraph{Coverage} In the range $n < 300$, the coverage of 
Theorem~\ref{t41} is as follows:
\begin{eqnarray*}
 35 =  7 \times  5: && 
  \uu{35}=\cy{11}{3}  \times \cy{34}{2} \times \cy{22}{4}      \\
 77 =  7 \times 11: && 
  \uu{77}=\cy{67}{3}  \times \cy{76}{2} \times \cy{8}{10}      \\
 95 = 19 \times  5: &&
  \uu{95}=\cy{16}{9}  \times \cy{94}{2} \times \cy{77}{4}      \\
119 =  7 \times 17: && 
  \uu{119}=\cy{18}{3}  \times \cy{118}{2} \times \cy{99}{16}   \\
155 = 31 \times  5: && 
  \uu{155}=\cy{121}{15}\times \cy{154}{2} \times \cy{32}{4}    \\
161 =  7 \times 23: && 
  \uu{161}=\cy{116}{3}  \times \cy{160}{2} \times \cy{43}{22}  \\
203 =  7 \times 29: && 
  \uu{203}=\cy{88}{3}  \times \cy{202}{2} \times \cy{113}{28}  \\
209 = 19 \times 11: && 
  \uu{209}=\cy{111}{9} \times \cy{208}{2} \times \cy{96}{10}   \\
215 = 43 \times  5: &&
  \uu{215}=\cy{126}{21}\times \cy{214}{2} \times \cy{87}{4}  
\end{eqnarray*}
The case $287=7\times41$ fails, since $\ord_{41}({-3}) = 8$.
In the range $q < 300$, the value $q = 41$ is the only prime $q$ with
$q \equiv 2$ \mod{3} and $\ord_q(-3) \neq q-1$.   However, in the range
$p < 300$, there are four primes $p$ with $p \equiv 7$ \mod{12} and
$\ord_p(-3) \neq (p-1)/2$, namely $p = 67$, 103, 151 and 271.
\bigskip

As is hinted in \S8 of \cite{cap}, many special cases arise when
we come to consider composite values $n = pq$
where $p$ and $q$ are distinct primes satisfying
$p \equiv q \equiv 1$ \mod{6}, with \textup{gcd}$(p-1,\,q-1) = 6$.
Accordingly, we do not offer theorems to cover these cases.
Instead, for the range $n < 1000$, we use Table~1 to list 
the instances in which we have
\begin{eqnarray}
\un = \cy{2x+3}{m} \times \cy{x+1}{3} \times \cy{-1}{2} 
      =   \cy{-2x-3}{m} \times \cy{-x-2}{3} \times \cy{-1}{2} 
\end{eqnarray}
where $m = \phi_n / 6$.   The following values of $n$ are \textbf{not} 
covered :
$259 = 7 \times 37$,\ \ $427 = 7 \times 61$,\ \
$511 = 7 \times 73$ and $973 = 7 \times 139$.

\begin{table}[p]
TABLE 1\\ \\
Decompositions (9) for $\un$ where $n = pq$ as specified in the text 
\begin{center}
$
\begin{array}
   {@{\hspace{1mm}}r@{\hspace{1mm}}c@{}c@{\hspace{2mm}}l@{\hspace{1mm}}l}
\hline
 91&=& 7 \times 13: &
   \uu{91} &=\cy{33}{12} \times \cy{16}{3} \times \cy{90}{2}
            =\cy{58}{12} \times \cy{74}{3} \times \cy{90}{2}    \\ 
133&=& 7 \times 19: & 
   \uu{133}&=\cy{61}{18} \times \cy{30}{3} \times \cy{132}{2}
            =\cy{72}{18} \times \cy{102}{3} \times \cy{132}{2}  \\
217&=& 7 \times 31: &
   \uu{217}&=\cy{135}{30} \times \cy{67}{3} \times \cy{216}{2}
            =\cy{82}{30} \times \cy{149}{3} \times \cy{216}{2}  \\
247&=&13 \times 19: & 
   \uu{247}&=\cy{137}{36} \times \cy{68}{3} \times \cy{246}{2}
            =\cy{110}{36} \times \cy{178}{3} \times \cy{246}{2} \\
   & &    &&=\cy{175}{36} \times \cy{87}{3} \times \cy{246}{2}
            =\cy{72}{36} \times \cy{159}{3} \times \cy{246}{2}  \\
301&=& 7 \times 43: &
   \uu{301}&=\cy{271}{42} \times \cy{135}{3} \times \cy{300}{2}
            =\cy{30}{42} \times \cy{165}{3} \times \cy{300}{2}  \\
403&=&13 \times 31: &
   \uu{403}&=\cy{228}{60} \times \cy{315}{3} \times \cy{402}{2}
            =\cy{175}{60} \times \cy{87}{3} \times \cy{402}{2}  \\
469&=& 7 \times 67: &
   \uu{469}&=\cy{142}{66} \times \cy{305}{3} \times \cy{468}{2}
            =\cy{327}{66} \times \cy{163}{3} \times \cy{468}{2} \\
553&=& 7 \times 79: &
   \uu{553}&=\cy{205}{78} \times \cy{102}{3} \times \cy{552}{2}
            =\cy{348}{78} \times \cy{450}{3} \times \cy{552}{2} \\
559&=&13 \times 43: &
   \uu{559}&=\cy{202}{84} \times \cy{380}{3} \times \cy{558}{2}
            =\cy{357}{84} \times \cy{178}{3} \times \cy{558}{2} \\
589&=&19 \times 31: &
   \uu{589}&=\cy{547}{90} \times \cy{273}{3} \times \cy{588}{2}
            =\cy{42}{90} \times \cy{315}{3} \times \cy{588}{2}  \\
679&=& 7 \times 97: &
   \uu{679}&=\cy{26}{96} \times \cy{352}{3} \times \cy{678}{2}
            =\cy{653}{96} \times \cy{326}{3} \times \cy{678}{2} \\
721&=& 7 \times 103:&
   \uu{721}&=\cy{422}{102} \times \cy{571}{3} \times \cy{720}{2}
            =\cy{299}{102} \times \cy{149}{3} \times \cy{720}{2}\\
763&=& 7 \times 109:&
   \uu{763}&=\cy{236}{108} \times \cy{499}{3} \times \cy{762}{2}
            =\cy{527}{108} \times \cy{263}{3} \times \cy{762}{2}\\
   & &    &&=\cy{345}{108} \times \cy{172}{3} \times \cy{762}{2}
            =\cy{418}{108} \times \cy{590}{3} \times \cy{762}{2}\\
817&=&19 \times 43:&
   \uu{817}&=\cy{357}{126} \times \cy{178}{3} \times \cy{816}{2}
            =\cy{460}{126} \times \cy{638}{3} \times \cy{816}{2}\\
871&=&13 \times 67:&
   \uu{871}&=\cy{59}{132} \times \cy{29}{3} \times \cy{870}{2}
            =\cy{812}{132} \times \cy{841}{3} \times \cy{870}{2}\\
   & &    &&=\cy{410}{132} \times \cy{640}{3} \times \cy{870}{2}
            =\cy{461}{132} \times \cy{230}{3} \times \cy{870}{2}\\
889&=& 7 \times127:&
   \uu{889}&=\cy{674}{126} \times \cy{781}{3} \times \cy{888}{2}
            =\cy{215}{126} \times \cy{107}{3} \times \cy{888}{2}\\
\hline
\end{array}
$
\end{center}
\end{table}

Now let $n = pq$ where $p$ and $q$ are distinct primes satisfying
$p \equiv q \equiv 5$ \mod{8}, $q > 5$ and gcd$(p-1,q-1) = 4$.
For the range $n < 1000$, Table~2 lists 3AP decompositions of~$\un$ 
that are lifts from weak 3AP decompositions of~$\uu{q}$; 
an asterisk marks a generator lifted from~1~\mod{q}.
Where the weak 3AP decomposition has a generator of order 4, 
we classify the lifted 3AP decompos\-itions into three types: 
if the generator lifted from~1 can be placed first, we
have type~A when the order of the middle generator is 4, and type~C
when the order of the last generator is 4, whereas type~B has the 
generator lifted from~1 in the middle.

\begin{table}[p]
TABLE 2\\
Some lifts from weak 3AP decompositions of $\uu{q}$, as specified
at the end of \S4.3
\renewcommand{\arraystretch}{1.1}
\begin{center}
$
\begin{array}{rcc}
  n = p \times  q & \textup{3AP decomposition of } \un &  \textup{Type}     \\
\hline
 65 = 5 \times 13 & \cy{27*}{4} \times \cy{44}{4} \times \cy{61}{3}     & A \\
                  & \cy{53*}{4} \times \cy{57}{4} \times \cy{61}{3}     & A \\
                  & \cy{53*}{4} \times \cy{16}{3} \times \cy{44}{4}     & C \\
145 = 5 \times 29 & \cy{88*}{4} \times \cy{12}{4} \times \cy{81}{7}     & A \\
                  & \cy{117*}{4} \times \cy{99}{4} \times \cy{81}{7}    & A \\
185 = 5 \times 37 & \cy{43}{4} \times \cy{112*}{4} \times \cy{181}{9}   & B \\
                  & \cy{38*}{4} \times \cy{16}{9} \times \cy{179}{4}    & C \\
265 = 5 \times 53 & \cy{213*}{4} \times \cy{201}{13} \times \cy{189}{4} & C \\
305 = 5 \times 61 & \cy{123*}{4} \times \cy{56}{15} \times \cy{294}{4}  & C \\
                  & \cy{62*}{4} \times \cy{24}{20} \times \cy{291}{3}   & - \\
                  & \cy{273}{12} \times \cy{62*}{4} \times \cy{156}{5}  & - \\
377 = 13\times 29 & \cy{262*}{12} \times \cy{99}{4} \times \cy{313}{7}  & A \\
377 = 29\times 13 & \cy{287*}{28} \times \cy{57}{4} \times \cy{203}{3}  & A \\
                  & \cy{14*}{28} \times \cy{146}{3} \times \cy{278}{4}  & C \\
                  & \cy{222*}{28} \times \cy{146}{3} \times \cy{70}{4}  & C \\
                  & \cy{235*}{28} \times \cy{146}{3} \times \cy{57}{4}  & C \\
505 = 5 \times 101& \cy{102*}{4} \times \cy{394}{4} \times \cy{181}{25} & A \\
                  & \cy{102*}{4} \times \cy{414}{4} \times \cy{221}{25} & A \\
                  & \cy{203*}{4} \times \cy{192}{4} \times \cy{181}{25} & A \\
                  & \cy{203*}{4} \times \cy{212}{4} \times \cy{221}{25} & A \\
                  & \cy{203*}{4} \times \cy{56}{25} \times \cy{414}{4}  & C \\
545 = 5 \times 109& \cy{33}{4} \times \cy{437*}{4} \times \cy{296}{27}  & B \\
                  & \cy{403}{4} \times \cy{437*}{4} \times \cy{471}{27} & B \\
689 = 13\times 53 & \cy{319*}{12} \times \cy{625}{13} \times \cy{242}{4}& C \\
689 = 53\times 13 & \cy{209*}{52} \times \cy{317}{4} \times \cy{425}{3} & A \\
                  & \cy{469*}{52} \times \cy{447}{4} \times \cy{425}{3} & A \\
                  & \cy{456*}{52} \times \cy{107}{3} \times \cy{447}{4} & C \\
                  & \cy{586*}{52} \times \cy{107}{3} \times \cy{317}{4} & C \\
745 = 5 \times 149& \cy{193}{4} \times \cy{597*}{4} \times \cy{256}{37} & B \\
                  & \cy{403}{4} \times \cy{597*}{4} \times \cy{46}{37}  & B \\
785 = 5 \times 157& \cy{158*}{4} \times \cy{757}{4} \times \cy{571}{39} & A \\
                  & \cy{472*}{4} \times \cy{129}{4} \times \cy{571}{39} & A \\
                  & \cy{443}{4} \times \cy{472*}{4} \times \cy{501}{39} & B \\
                  & \cy{158*}{4} \times \cy{207}{12} \times \cy{256}{13}& - \\
                  & \cy{158*}{4} \times \cy{326}{3} \times \cy{494}{52} & - \\
865 = 5 \times 173& \cy{693*}{4} \times \cy{566}{43} \times \cy{439}{4} & - \\
905 = 5 \times 181& \cy{363*}{4} \times \cy{316}{5} \times \cy{269}{36} & - \\
985 = 5 \times 197& \cy{183}{4} \times \cy{592*}{4} \times \cy{16}{49}  & B \\
\hline
\end{array}
$
\end{center}
\renewcommand{\arraystretch}{1.4}
\end{table}

\subsection{Lifts from $n = kp$ to $n = kp^2$}

For the range $n < 1000$, details of the 3AP decompositions (3APDs) for 
values of the form $n = kp^2$ ($k$ and $p$ distinct odd primes, $k > 3$,
$p > 3$) are as in Table~3.
With the given restrictions on $k$ and $p$,
just one value of the form $n = kp^3$ lies in the range $n < 1000$,
namely $n = 875$, and it has precisely six 3AP decompositions.
Each of these is obtained by further lifting one of the 3AP decompositions
for $n = 175$.   In this further lifting, the orders that are not
multiples of 5 are unchanged, but the orders that are multiples of 5
become multiples of $5^2$.

\begin{table}[p]
TABLE 3\\ \\
An enumeration of decompositions for $n = kp^2$\\
\begin{center}
\renewcommand{\arraystretch}{1.1}
$
\begin{array}{lr@{\hspace{5mm}}rrr}
        && \multicolumn{3}{c}{\overbrace{\hspace{60mm}}}           \\
\ \ \ n = kp^2 & \#\textup{ 3APDs}  
         & \multicolumn{2}{c}{\#\textup{ lifts from }\uu{kp}}
         & \#\textup{ other}                                       \\
\cline{3-4}
        && \multicolumn{1}{c}{\textup{from}}
         & \multicolumn{1}{c}{\textup{from}} &                     \\
        && \textup{strong 3APDs}
         & \textup{weak 3APDs} &                                    \\
\hline
175 = 7 \cdot 5^2    &   6 &   3 &   3 &   0 \\
245 = 5 \cdot 7^2    &   6 &   3 &   3 &   0 \\
275 = 11 \cdot 5^2   &  68 &   0 &   8 &  60 \\
325 = 13 \cdot 5^2   &  20 &  12 &   8 &   0 \\
425 = 17 \cdot 5^2   &   8 &   0 &   8 &   0 \\
475 = 19 \cdot 5^2   &   6 &   3 &   3 &   0 \\
539 = 11 \cdot 7^2   &  12 &   9 &   3 &   0 \\
575 = 23 \cdot 5^2   &   2 &   0 &   2 &   0 \\
605 = 5 \cdot 11^2   &   0 &   0 &   0\rlap{*} & 0 \\
637 = 13 \cdot 7^2   & 126 & 108 &  18 &   0 \\
725 = 29 \cdot 5^2   &  30 &  18 &  12 &   0 \\
775 = 31 \cdot 5^2   & 188 &  32 &  24 & 132 \\
845 = 5 \cdot 13^2   &  20 &  12 &   8 &   0 \\
847 = 7 \cdot 11^2   &   0 &   0\rlap{*} &   0\rlap{*} &   0 \\
925 = 37 \cdot 5^2   &  10 &   6 &   4 &   0 \\
931 = 19 \cdot 7^2   & 182 & 156 &  26 &   0 \\
\hline\\
\end{array}
$\\ 
\renewcommand{\arraystretch}{1.4}
{*} 3APDs of $\uu{kp}$ exist, but the special lifts of the generators 
    are in AP
\end{center}
\end{table}

As Table~3 indicates, some of the 3AP decompositions for $n = 275$ 
and $775$ are not lifts, these two $n$-values being distinctive in that
they have $ p\ |\ (k-1) $\,.   How do these exceptional 3AP decompositions
arise?   One of them is given by 
\[ \uu{275} = \cy{16}{5} \times \cy{24}{10} \times \cy{32}{4}\,.         \]
This is related to the decomposition
\[ \uu{55}  = \phantom{\cy{116}{5} \times} \cy{24}{10} \times \cy{32}{4} \]
and to the fact that, within $\uu{55}$, we have 
$\cy{16}{5} \subset \cy{24}{10}$.   The further 3AP decomposition
\[ \uu{275} = \cy{181}{5} \times \cy{244}{10} \times \cy{32}{4}          \]
arises in the same way, as 181 is a lift of 16, and 244 is a lift of 24.

Less straightforward situtations exist too.   Consider, for example, the
3AP decomposition
\[ \uu{775} = \cy{32}{4} \times \cy{54}{10} \times \cy{76}{15}           \]
with $n = 31 \cdot 5^2$.   If we try lifting to this from 
$n = 31 \cdot 5 = 155$,
we find that 32, 54 and 76 also have orders 4, 10 and 15 \mod{155}, and that
\begin{eqnarray*}
 \uu{155} & = & \cy{32}{4} \times \cy{54}{10} \times \cy{76^5}{3}        \\
          & = & \cy{32}{4} \times \cy{54^5}{2} \times \cy{76}{15}\,.
\end{eqnarray*}
Analogous to this, we can rewrite an example from the previous paragraph
in the weak form
\[ \uu{55}  = \cy{16^5}{1} \times \cy{24}{10} \times \cy{32}{4}\,.       \]
These examples suggest an amusing generalisation of 3AP decompositions to
decompositions of the form
\[ \un = \cy{x^h}{} \times \cy{(x+k)^i}{} \times \cy{(x+2k)^j}{}         \]
where $h,\,i,\,j \ge 1\,$, but we do not pursue this idea further here.

\section{A class of weak 3AP decompositions}
\label{s5}

A noteworthy class of weak 3AP decompositions arises for primes $n$ 
satisfying $n \equiv 1$ \mod{6p} where $p$ is an odd prime, $p > 3$.
In each of these decompositions, one of the generators has order 6 and
another has order~$p$.   For $n < 300$, such decompositions are as follows:
\begin{eqnarray*}
n =  43  && \cy{1}{1}\times\cy{4}{7}\times\cy{7}{6}         \\
n =  67  && \cy{1}{1}\times\cy{30}{6}\times\cy{59}{11}      \\
n =  79  && \cy{1}{1}\times\cy{52}{13}\times\cy{24}{6}      \\
n = 103  && \cy{1}{1}\times\cy{47}{6}\times\cy{93}{17}      \\
n = 139  && \cy{97}{6}\times\cy{1}{1}\times\cy{44}{23}      \\
n = 223  & \left\{ \begin{array}{@{}c@{}}  \\ \\ \end{array} \right.
         & \begin{array}{@{}l@{}}
            \cy{1}{1}\times\cy{132}{37}\times\cy{40}{6}  \\
            \cy{184}{6}\times\cy{1}{1}\times\cy{41}{37}  \\
            \end{array}                                     \\  
n = 283  && \cy{45}{6}\times\cy{1}{1}\times\cy{240}{47}
\end{eqnarray*}
Where $\un = \cy{a}{6}\times\cy{1}{1}\times\cy{c}{p}$, we have
$a \equiv (c-1)^{-1}$ \mod{n}, so that
$\un = $ \mbox{$\cy{c-1}{6}\times\cy{c}{p}$}, a situation discussed in
\cite[\S8.2]{cap}.   The above weak 3AP decomposition for $n = 67$
has the \emph{orders} of the generators in AP.

\section{Finite fields}
\label{s6}

Finite fields of non-prime order can have 3AP decompositions. 
Clearly this is impossible in fields of characteristic~$2$,
which contain no 3-term arithmetic progressions. 

\paragraph{Example} 
Using GAP~\cite{GAP}, we found the following 3AP decompositions 
of small finite fields $\GF(q)$.
In this list, $\zeta$ denotes the primitive root 
(denoted by \verb+Z(q)+ in GAP notation) in the field $\GF(q)$.
It is a root of the appropriate \emph{Conway polynomial} \cite{Conway};
the relevant Conway polynomials are as follows:
\begin{eqnarray*}
q=11^2: && x^2+7x+2 \\
q=11^3: && x^3+2x+9 \\
q=19^2: && x^2-x+2 \\
q=19^3: && x^3+4x-2 \\
q=23^2: && x^2-2x+5 \\
q=29^2: && x^2-5x+2
\end{eqnarray*}
We give only
one decomposition for each possible list of orders of the factors:
\begin{eqnarray*}
\GF(11^2)^\times 
&=& \cy{\zeta^{72}}{5} \times \cy{\zeta^{15}}{8} \times \cy{\zeta^{80}}{3} \\
\GF(11^3)^\times 
&=& \cy{\zeta^{570}}{7} \times \cy{\zeta^{532}}{5} \times \cy{\zeta^{595}}{38} \\
&=& \cy{\zeta^{665}}{2} \times \cy{\zeta^{1008}}{95} \times \cy{\zeta^{570}}{7} \\
\GF(19^2)^\times 
&=& \cy{\zeta^{144}}{5} \times \cy{\zeta^{320}}{9} \times \cy{\zeta^{135}}{8} \\
\GF(19^3)^\times 
&=& \cy{\zeta^{3429}}{2} \times \cy{\zeta^{2970}}{127} \times \cy{\zeta^{5588}}{27} \\
\GF(23^2)^\times 
&=& \cy{\zeta^{176}}{3} \times \cy{\zeta^{192}}{11} \times \cy{\zeta^{429}}{16} \\
\GF(29^2)^\times 
&=& \cy{\zeta^{280}}{3} \times \cy{\zeta^{720}}{7} \times \cy{\zeta^{609}}{40}\\
&=& \cy{\zeta^{120}}{7} \times \cy{\zeta^{504}}{5} \times \cy{\zeta^{385}}{24} 
\end{eqnarray*}

Can we have a 3AP decomposition of
$\GF(q)^\times$ in which two of the generators have orders $2$ and $3$\,? 
As earlier, such a decomposition requires that $(q-1)/6$ is co-prime to $6$,
so that $q\equiv7$ or $31$ \mod{36}. But this implies that, if $q=p^n$ with
$p$ prime, then $p\equiv7$ or $31$ \mod{36} (since $7$ and $31$ are 
\PLRs of 36 \cite{cap}, and each is the fifth power of the other). Then 
elements of orders $2$ and $3$ lie in the prime subfield, and hence 
so does the whole AP. So there are no such decompositions 
other than those of $\un$ for $n$ prime discussed in \S\ref{s2}.

A similar argument shows that a 3AP decomposition of $\GF(11^3)^\times$ 
into \mbox{factors} of orders $2$, $5$ and $133$ is impossible.

\section{Decompositions with more than three factors}
\label{s7}

As indicated above, we can define 4AP decompositions analogously 
to 3AP decompositions.

A computer program has shown that
no examples of strong 4AP decompositions of $\un$ exist
for prime values of $n$ up to $10000$.  The smallest composite $n$ 
for which strong 4AP decompositions of $\un$ exist is even: 
\begin{eqnarray*}
\uu{104} & = & 
\cy{31}{4}\times\cy{81}{3}\times\cy{27}{2}\times\cy{77}{2} \\
                 & = & 
\cy{77}{2}\times\cy{79}{2}\times\cy{81}{3}\times\cy{83}{4}\,.
\end{eqnarray*}
The smallest weak $4$AP decomposition of $\un$ with prime $n$ is
\[\uu{3613} = 
\cy{3528}{4} \times \cy{1148}{129} \times \cy{2381}{7} \times\cy{1}{1}\,.\]
We have no examples with larger numbers of generators that are in AP.

\paragraph{Note} The computations reported in this paper were performed
using GAP~\cite{GAP}, and a package of GAP functions written by the first
author for computations in the groups \un, available from~\cite{plrfns}.
Further documentation of these functions can be found in~\cite{cap}.

\end{document}